\documentclass[12pt,english]{article}
\usepackage[T1]{fontenc}
\usepackage[latin1]{inputenc}
\usepackage[a4paper]{geometry}
\usepackage{babel}
\usepackage{prettyref}
\usepackage{mathrsfs}
\usepackage{amsthm}
\usepackage{amsmath}
\usepackage{amssymb}
\usepackage{esint}
\usepackage[numbers,compress]{natbib}
\usepackage[unicode=true,pdfusetitle,
 bookmarks=true,bookmarksnumbered=false,bookmarksopen=false,
 breaklinks=false,pdfborder={0 0 1},backref=false,colorlinks=false]
 {hyperref}

\makeatletter
\theoremstyle{plain}
\newtheorem{thm}{\protect\theoremname}
  \theoremstyle{plain}
  \newtheorem{lem}[thm]{\protect\lemmaname}
  \theoremstyle{plain}
  \newtheorem*{thm*}{\protect\theoremname}
  \theoremstyle{remark}
  \newtheorem*{rem*}{\protect\remarkname}

\@ifundefined{date}{}{\date{}}
\usepackage[colorinlistoftodos]{todonotes}

\makeatother

  \providecommand{\lemmaname}{Lemma}
  \providecommand{\remarkname}{Remark}
  \providecommand{\theoremname}{Theorem}
\providecommand{\theoremname}{Theorem}

\begin{document}

\title{\textbf{Structural stability for a thermal convection model with
temperature-dependent solubility}}

\author{M. Ciarletta%
\thanks{Dipartimento di Ingegneria Industriale, Università di Salerno, Italy,
\texttt{mciarletta@unisa.it}%
}\and{}B. Straughan%
\thanks{Department of Mathematical Sciences, University of Durham, UK, \texttt{brian.straughan@durham.ac.uk}%
}\and{}V. Tibullo%
\thanks{Corresponding author, Dipartimento di Ingegneria dell'Informazione,
Ingegneria Elettrica e Matematica Applicata, Università di Salerno,
Italy, \texttt{vtibullo@unisa.it}%
}}
\maketitle
\begin{abstract}
We study a problem involving thermosolutal convection in a fluid when
the solute concentration is subject to a chemical reaction in which
the solubility of the dissolved component is a function of temperature.
When the spatial domain is a bounded one in $\mathbb{R}^{2}$ we show
that the solution depends continuously on the reaction rate using
true \emph{a priori} bounds for the solution when the chemical equilibrium
function is an arbitrary function of temperature.

Link to publisher: \url{http://dx.doi.org/10.1016/j.nonrwa.2014.07.012}
\end{abstract}

\section{Introduction}

The problem of double diffusion convection in a horizontal layer of
fluid with simultaneous chemical reaction is well studied, cf. the
account in \citet{Straughan:2004}, pp. 225-237. Indeed the governing
equations for such a chemical reaction problem are derived by \citet{MorroStraughan:1990},
employing ideas of continuum thermodynamics. It is also worth drawing
attention to the fact that similar equations arise in fluid phase
change problems where continuum thermodynamic theories are again employed
and analysed by \citet{Berti:2014}, \citet{Berti:2013}, \citet{Bonetti:2009},
\citet{Bonetti:2011}, \citet{Fabrizio:2012} and \citet{Fabrizio:2006,Fabrizio:2008,Fabrizio:2011}.
A particular reaction where the solute concentration is subject to
a chemical reaction in which the solubility of the dissolved component
is a linear function of the temperature has been analysed recently
in a porous medium context by \citet{PritchardRichardson:2007}, by
\citet{WangTan:2009}, and by \citet{MalashettyBiradar:2011}. The
object of this paper is to analyse the effect of such a class of reaction
terms in a fluid. The dependence on temperature is here taken to be
arbitrary and not only linear. We allow the fluid to occupy a bounded
domain $\Omega\subset\mathbb{R}^{2}$ with boundary $\Gamma$ sufficiently
smooth to allow application of the divergence theorem. We are particularly
interested to investigate the continuous dependence of the solution
on the reaction rate. It is very reasonable to expect that in such
a mathematical model one can achieve an appropriate type of continuous
dependence and this kind of stability problem belongs to the important
class of structural stability questions.

Structural stability, or continuous dependence on the model itself,
is a concept al least as important as the classical idea of stability
which involves continuous dependence on the initial data, as explained
in some detail by \citet{HirschSmale:1974}. Structural stability
was studied in elasticity by \citet{KnopsPayne:1969}, then further
advanced in a variety of continuum mechanical contexts by \citet{Payne:1987,Payne:1987b,Payne:1989}
and by \citet{KnopsPayne:1988}. Since then structural stability studies
in continuum mechanics have proved very popular as witnessed by the
works of \citet{Aulisa:2009,Aulisa:2011}, \citet{Celebi:2006}, \citet{Chirita:2006},
\citet{Ciarletta:2011,Ciarletta:2012}, \citet{HoangIbragimov:2011,HoangIbragimov:2012},
\citet{Hoang:2013}, \citet{KalantarovZelik:2012}, \citet{Kandem:2011},
\citet{KangPark:2013}, \citet{Kelliher:2011}, \citet{Li:2012},
\citet{Liu:2009,Liu:2012}, \citet{Liu:2010,Liu:2010b}, \citet{OuyangYang:2009},
\citet{Passarella:2014}, \citet{Ugurlu:2008}, \citet{You:2012}.

\section{Fundamental equations}

The fundamental model we study is based upon the equations of balance
of momentum, balance of mass, conservation of energy, and conservation
of salt concentration, adopting a Boussinesq approximation in the
body force term in the momentum equation. Thus, let $v_{i}(\mathbf{x},t)$,
$p(\mathbf{x},t)$, $T(\mathbf{x},t)$ and $C(\mathbf{x},t)$ denote
velocity, pressure, temperature and salt concentration, where $\mathbf{x}\in\Omega$,
$t$ denote time, with $0<t<\mathscr{T}$, for some $\mathscr{T}<\infty$.
Then, the equations of momentum, mass, energy, and salt concentration
are taken to be
\begin{equation}
\begin{alignedat}{1}\frac{\partial v_{i}}{\partial t}+v_{j}\frac{\partial v_{i}}{\partial x_{j}} & =-\frac{\partial p}{\partial x_{i}}+\Delta v_{i}+g_{i}T-h_{i}C,\\
\frac{\partial v_{i}}{\partial x_{i}} & =0,\\
\frac{\partial T}{\partial t}+v_{i}\frac{\partial T}{\partial x_{i}} & =\Delta T,\\
a\frac{\partial C}{\partial t}+bv_{i}\frac{\partial C}{\partial x_{i}} & =\Delta C+Lf(T)-KC,
\end{alignedat}
\label{eq:fundamental}
\end{equation}
where $g_{i}$, $h_{i}$ represent gravity vectors and without loss
of generality we assume that
\[
|\mathbf{g}|,\,|\mathbf{h}|\leq1.
\]
Standard indicial notation is used throughout with a repeated index
denoting summation over $1$ and $2$, and $a$, $b$, $L$ and $K$
are positive constants. Equations \eqref{eq:fundamental} follow in
practice by employing a Boussinesq approximation which accounts for
variable $C$ allowing the incompressibility condition to hold, cf.
\citet{Fife:2000}, pp. 72-74. 

The function $f$ is at least $C^{1}$ and the term $Lf(T)$ is analogous
to the chemical equilibrium term, $C_{eq}$, in \citet{PritchardRichardson:2007},
\citet{WangTan:2009}, and \citet{MalashettyBiradar:2011}, although
all of these writers assume $C_{eq}(T)$ is a linear function of $T$.
The terms $Lf(T)-KC$ in equation \eqref{eq:fundamental}$_{4}$ correspond
to the mass supply term $m_{\alpha}$ in \citet{MorroStraughan:1990}.
The justification for this is as given in \citet{PritchardRichardson:2007}
who write $Lf(T)-KC=K(C_{eq}-C)$, $C_{eq}$ being a chemical equilibrium
term. The logic is that in chemical equilibrium the chemical reaction
arises solely due to the term $K(C_{eq}(T)-C)$.

In general \eqref{eq:fundamental}$_{4}$ holds with $a=b=1$. However,
we allow for different $a$ and $b$ since these coefficients will
change under a rescaling, i.e. under a different non-dimensionalization. 

Equations \eqref{eq:fundamental} hold in the domain $\Omega\times(0,\mathscr{T})$,
together with the initial conditions,
\begin{equation}
\begin{alignedat}{1}v_{i}(\mathbf{x},0) & =v_{i}^{0}(\mathbf{x}),\\
T(\mathbf{x},0) & =T_{0}(\mathbf{x}),\\
C(\mathbf{x},0) & =C_{0}(\mathbf{x}),
\end{alignedat}
\label{eq:initial-conditions}
\end{equation}
for $\mathbf{x}\in\Omega$, and the boundary conditions
\begin{equation}
\begin{alignedat}{1}v_{i}(\mathbf{x},t) & =0,\\
T(\mathbf{x},t) & =g(\mathbf{x},t),\\
C(\mathbf{x},t) & =h(\mathbf{x},t),
\end{alignedat}
\label{eq:boundary-conditions}
\end{equation}
$\mathbf{x}\in\Gamma$, $t\in[0,\mathscr{T})$.

Since we are interested in studying continuous dependence on the reaction
rates $L$ and $K$, we let $(u_{i},p_{1},T_{1},C_{1})$ and $(v_{i},p_{2},T_{2},C_{2})$
be two solutions to \eqref{eq:fundamental}-\eqref{eq:boundary-conditions}
for the same initial and boundary data, but for different reaction
coefficients $(L_{1},K_{1})$ and $(L_{2},K_{2})$. To progress we
now introduce the difference variables $(w_{i},\pi,\theta,\phi)$
and $l$ and $k$ by
\begin{equation}
\begin{aligned}w_{i} & =u_{i}-v_{i}, & \pi & =p_{1}-p_{2}, & \theta & =T_{1}-T_{2},\\
\phi & =C_{1}-C_{2}, & l & =L_{1}-L_{2}, & k & =K_{1}-K_{2}.
\end{aligned}
\label{eq:difference-variables}
\end{equation}

Thus, from \eqref{eq:fundamental}-\eqref{eq:difference-variables}
we may determine the boundary-initial value problem for the difference
variables as
\begin{equation}
\begin{alignedat}{1}\frac{\partial w_{i}}{\partial t}+w_{j}\frac{\partial u_{i}}{\partial x_{j}}+v_{j}\frac{\partial w_{i}}{\partial x_{j}} & =-\frac{\partial\pi}{\partial x_{i}}+\Delta w_{i}+g_{i}\theta-h_{i}\phi,\\
\frac{\partial w_{i}}{\partial x_{i}} & =0,\\
\frac{\partial\theta}{\partial t}+w_{i}\frac{\partial T^{1}}{\partial x_{i}}+v_{i}\frac{\partial\theta}{\partial x_{i}} & =\Delta\theta,\\
a\frac{\partial\phi}{\partial t}+b\left(w_{i}\frac{\partial C^{1}}{\partial x_{i}}+v_{i}\frac{\partial\phi}{\partial x_{i}}\right) & =\Delta\phi+L_{1}[f(T^{1})-f(T^{2})]+lf(T^{2})-K_{1}\phi-kC_{2},
\end{alignedat}
\label{eq:difference-equations}
\end{equation}
in $\Omega\times(0,\mathscr{T})$, with
\begin{equation}
w_{i}(\mathbf{x},0)=0,\quad\theta(\mathbf{x},0)=0,\quad\phi(\mathbf{x},0)=0,
\end{equation}
$\mathbf{x}\in\Omega$, together with
\begin{equation}
w_{i}(\mathbf{x},t)=0,\quad\theta(\mathbf{x},t)=0,\quad\phi(\mathbf{x},t)=0,\label{eq:difference-boundary-conditions}
\end{equation}
for $\mathbf{x}\in\Gamma$, $t\in[0,\mathscr{T})$. We write $T_{1}\equiv T^{1}$
when the occasion needs, with similar notation for $T_{2}$, $C_{1}$
and $C_{2}$.

We wish to derive a continuous dependence estimate for a suitable
measure of $w_{i}$, $\theta$, $\phi$ in terms of $l$ and $k$.
However, we require this estimate to be truly \emph{a priori} in the
sense that the coefficients which appear associated to $l$ and $k$
involve only data. Thus, before we can achieve our continuous dependence
result we need to derive some \emph{a priori} bounds for the solution
to \eqref{eq:fundamental}-\eqref{eq:boundary-conditions}.

\section{\emph{A priori} estimates}

We denote by $(\cdot,\cdot)$ and $\left\Vert \cdot\right\Vert $
the inner product and norm on $L^{2}(\Omega)$ and further let $\left\Vert \cdot\right\Vert _{p}$
be the norm on $L^{p}(\Omega)$ with $\left\Vert \cdot\right\Vert _{\infty}$
being the $L^{\infty}(\Omega)$ norm.

Define the quantity $T_{m}$ by
\begin{equation}
T_{m}=\max\{\left\Vert T_{0}\right\Vert _{\infty},\sup_{[0,\mathscr{T}]}\left\Vert g\right\Vert _{\infty}\}.
\end{equation}
\citet{Payne:2001} show how we may employ the function
\[
\tilde{\psi}=[T-T_{m}]^{+}=\sup(T-T_{m},0)
\]
to show from \eqref{eq:fundamental}-\eqref{eq:boundary-conditions}
that
\begin{equation}
\sup_{\Omega\times[0,\mathscr{T}]}|T(\mathbf{x},t)|\leq T_{m}.\label{eq:temperature-estimate}
\end{equation}
This \emph{a priori} bound for the temperature is very important in
what follows.

We need to establish \emph{a priori} bounds for certain norms of $\mathbf{v}$
and $C$ and to do this we recall two lemmas. The first arises from
a Rellich identity as used by \citet{PayneWeinberger:1958}, and is
given explicitly in \citet{PayneStraughan:1998jmpa}, inequality $(\text{A}10)$.
\begin{lem}
\label{lem:1}Let $\Phi$ be a harmonic function in $\Omega$ with
boundary values $Q$, i.e. $\Phi$ satisfies
\begin{equation}
\begin{alignedat}{2}\Delta\Phi & =0 & \qquad & \text{in }\Omega,\\
\Phi & =Q &  & \text{on \ensuremath{\Gamma}.}
\end{alignedat}
\label{eq:lemma-1-hyp}
\end{equation}
Then one may derive explicit constants $c_{1}$ and $c_{2}$ such
that
\begin{equation}
\left\Vert \nabla\Phi\right\Vert ^{2}+c_{1}\oint_{\Gamma}\left(\frac{\partial\Phi}{\partial n}\right)^{2}dA\leq c_{2}\oint_{\Gamma}|\nabla_{s}Q|^{2}dA,
\end{equation}
where $\nabla_{s}$ denotes the tangential derivative.
\end{lem}
\prettyref{lem:1} holds for a general domain $\Omega$ not just one
in $\mathbb{R}^{2}$. In the present case since we are in a two-dimensional
domain it is to be understood that the integral element $dA$ stands
for an integral along a curve. We use this notation consistently.
The second lemma is given as inequality $(\text{A}12)$ in \citet{PayneStraughan:1998jmpa}.
\begin{lem}
\label{lem:2}Let $\psi$ be the torsion function which satisfies
the boundary value problem
\begin{equation}
\begin{alignedat}{2}\Delta\psi & =-1 & \qquad & \text{in }\Omega,\\
\psi & =0 &  & \text{on \ensuremath{\Gamma}.}
\end{alignedat}
\end{equation}
Then by the maximum principle $\psi>0$ in $\Omega$, and for a function
$\Phi$ satisfying equations \eqref{eq:lemma-1-hyp} we have the inequality
\begin{equation}
2(\psi\nabla\Phi,\,\nabla\Phi)+\left\Vert \Phi\right\Vert ^{2}\leq\psi_{1}\oint_{\Gamma}Q^{2}dA,
\end{equation}
where 
\[
\psi_{1}=\max_{\Gamma}\left|\frac{\partial\psi}{\partial n}\right|.
\]

\end{lem}
We next take the inner product on $L^{2}(\Omega)$ of equation \eqref{eq:fundamental}$_{1}$
with $v_{i}$ and using the boundary conditions \eqref{eq:boundary-conditions},
integration by parts, and the Cauchy-Schwarz, arithmetic-geometric
mean and Poincaré inequalities, we obtain
\begin{equation}
\begin{alignedat}{1}\frac{d}{dt}\frac{1}{2}\left\Vert \mathbf{v}\right\Vert ^{2}+\left\Vert \nabla\mathbf{v}\right\Vert ^{2} & \leq\left\Vert T\right\Vert \left\Vert \mathbf{v}\right\Vert +\left\Vert C\right\Vert \left\Vert \mathbf{v}\right\Vert \\
 & \leq(\left\Vert T\right\Vert ^{2}+\left\Vert C\right\Vert ^{2})\frac{\alpha}{2}+\frac{1}{\alpha}\left\Vert \mathbf{v}\right\Vert ^{2}\\
 & \leq(\left\Vert T\right\Vert ^{2}+\left\Vert C\right\Vert ^{2})\frac{\alpha}{2}+\frac{1}{\alpha\lambda_{1}}\left\Vert \nabla\mathbf{v}\right\Vert ^{2},
\end{alignedat}
\label{eq:proof-14}
\end{equation}
where $\alpha>0$ is to be chosen and $\lambda_{1}$ is the first
eigenvalue in the membrane problem for $\Omega$. Pick now $\alpha=2/\lambda_{1}$
and integrate \eqref{eq:proof-14} to find
\begin{equation}
\left\Vert \mathbf{v}\right\Vert ^{2}+\int_{0}^{t}\left\Vert \nabla\mathbf{v}\right\Vert ^{2}ds\leq\left\Vert \mathbf{v}_{0}\right\Vert ^{2}+\frac{2}{\lambda_{1}}\int_{0}^{t}\left\Vert T\right\Vert ^{2}ds+\frac{2}{\lambda_{1}}\int_{0}^{t}\left\Vert C\right\Vert ^{2}ds.\label{eq:proof-15}
\end{equation}
Thanks to the estimate \eqref{eq:temperature-estimate} we may replace
the first two terms on the right of \eqref{eq:proof-15} by the data
term
\begin{equation}
d_{5}=\left\Vert \mathbf{v}_{0}\right\Vert ^{2}+\frac{2m\mathscr{T}T_{m}^{2}}{\lambda_{1}},\label{eq:d5-definition}
\end{equation}
where $m=m(\Omega)$ is the Lebesgue measure of $\Omega$.

Thus, from \eqref{eq:proof-15}, \eqref{eq:d5-definition} we may
establish
\begin{equation}
\left\Vert \mathbf{v}\right\Vert ^{2}+\int_{0}^{t}\left\Vert \nabla\mathbf{v}\right\Vert ^{2}ds\leq d_{5}+\frac{2}{\lambda_{1}}\int_{0}^{t}\left\Vert C\right\Vert ^{2}ds,\label{eq:proof-17}
\end{equation}
and from this with the aid of Poincaré's inequality we see that
\begin{equation}
\left\Vert \mathbf{v}\right\Vert ^{2}+\lambda_{1}\int_{0}^{t}\left\Vert \mathbf{v}\right\Vert ^{2}ds\leq d_{5}+\frac{2}{\lambda_{1}}\int_{0}^{t}\left\Vert C\right\Vert ^{2}ds.\label{eq:proof-18}
\end{equation}

From this point we introduce the function $H(\mathbf{x},t)$ by
\begin{equation}
\begin{alignedat}{2}\Delta H & =0 & \qquad & \text{in }\Omega,\\
H & =h &  & \text{on \ensuremath{\Gamma},}
\end{alignedat}
\label{eq:h-problem}
\end{equation}
where $h(\mathbf{x},t)$ is the boundary data function for $C$ as
given in equations \eqref{eq:boundary-conditions}.

Now, multiply equation \eqref{eq:fundamental}$_{4}$ by $C-H$ and
integrate to find
\begin{equation}
\begin{gathered}a\int_{0}^{t}\int_{\Omega}C_{,s}(C-H)dxds=-b\int_{0}^{t}\int_{\Omega}v_{i}C_{,i}(C-H)dxds\\
+\int_{0}^{t}\int_{\Omega}\Delta C(C-H)dxds+L\int_{0}^{t}\int_{\Omega}f(T)(C-H)dxds\\
-K\int_{0}^{t}\int_{\Omega}C(C-H)dxds.
\end{gathered}
\label{eq:proof-20}
\end{equation}
Denote the five terms in equation \eqref{eq:proof-20} by $I_{1}$--$I_{5}$
and we now develop these. By integration
\begin{equation}
\begin{gathered}I_{1}=\frac{a}{2}(\left\Vert C\right\Vert ^{2}-\left\Vert C_{0}\right\Vert ^{2})-a\int_{\Omega}HCdx+a\int_{\Omega}H_{0}C_{0}dx\\
+a\int_{0}^{t}\int_{\Omega}CH_{,s}dxds.
\end{gathered}
\label{eq:i1}
\end{equation}

Using the maximum principle $H$ may be bounded by $h_{m}$ where
\[
h_{m}=\max_{\Gamma\times[0,\mathscr{T}]}|h|.
\]
Then
\begin{equation}
\begin{gathered}I_{2}=+b\int_{0}^{t}\int_{\Omega}v_{i}C_{,i}Hdxds\leq bh_{m}\sqrt{\int_{0}^{t}\left\Vert \mathbf{v}\right\Vert ^{2}ds\,\int_{0}^{t}\left\Vert \nabla C\right\Vert ^{2}ds}\\
\leq bh_{m}\sqrt{\frac{d_{5}}{\lambda_{1}}+\frac{2}{\lambda_{1}^{2}}\int_{0}^{t}\left\Vert C\right\Vert ^{2}ds}\,\sqrt{\int_{0}^{t}\left\Vert \nabla C\right\Vert ^{2}ds},
\end{gathered}
\label{eq:i2}
\end{equation}
where in deriving \eqref{eq:i2} we have employed the Cauchy-Schwarz
inequality and estimate \eqref{eq:proof-18}.

For $I_{3}$ we integrate by parts and use \eqref{eq:h-problem} to
obtain
\[
\begin{alignedat}{1}I_{3} & =-\int_{0}^{t}\left\Vert \nabla C\right\Vert ^{2}ds+\int_{0}^{t}(\nabla C,\,\nabla H)ds\\
 & =-\int_{0}^{t}\left\Vert \nabla C\right\Vert ^{2}ds+\int_{0}^{t}\oint_{\Gamma}h\frac{\partial H}{\partial n}dAds.
\end{alignedat}
\]
We further employ the Cauchy-Schwarz inequality, and then \prettyref{lem:1}
to find
\begin{equation}
I_{3}\leq-\int_{0}^{t}\left\Vert \nabla C\right\Vert ^{2}ds+\sqrt{\frac{c_{2}}{c_{1}}\int_{0}^{t}\oint_{\Gamma}h^{2}dAds\,\int_{0}^{t}\oint_{\Gamma}|\nabla_{s}h|^{2}dAds}.
\end{equation}
To estimate $I_{4}$ we use the arithmetic-geometric mean inequality
with positive constants $\gamma_{1}$ and $\gamma_{2}$ to find
\[
I_{4}\leq\frac{L}{2}(\gamma_{1}^{-1}+\gamma_{2}^{-1})\int_{0}^{t}\left\Vert f\right\Vert ^{2}ds+\frac{L\gamma_{1}}{2}\int_{0}^{t}\left\Vert C\right\Vert ^{2}ds+\frac{L\gamma_{2}}{2}\int_{0}^{t}\left\Vert H\right\Vert ^{2}ds.
\]
Now, $f$ is known and $f=f(T)$ and so using bound \prettyref{eq:temperature-estimate}
we may bound $\int_{0}^{t}\left\Vert f\right\Vert ^{2}ds$ by data,
say $d_{4}.$ Further, employ \prettyref{lem:2} on the $\left\Vert H\right\Vert $
term and we then find
\begin{equation}
I_{4}\leq\frac{L}{2}(\gamma_{1}^{-1}+\gamma_{2}^{-1})d_{4}+\frac{L\gamma_{1}}{2}\int_{0}^{t}\left\Vert C\right\Vert ^{2}ds+\frac{L\gamma_{2}\psi_{1}}{2}\int_{0}^{t}\oint_{\Gamma}h^{2}dAds.
\end{equation}
Finally employing the arithmetic-geometric mean inequality with a
constant $\zeta_{3}>0$ together with \prettyref{lem:2} we obtain
\begin{equation}
I_{5}\leq-\left(K-\frac{K}{2\zeta_{3}}\right)\int_{0}^{t}\left\Vert C\right\Vert ^{2}ds+\frac{\psi_{1}K\zeta_{3}}{2}\int_{0}^{t}\oint_{\Gamma}h^{2}dAds.\label{eq:i5}
\end{equation}

We now group \eqref{eq:i1}-\eqref{eq:i5} together in equation \eqref{eq:proof-20}
and with further use of the arithmetic-geometric mean inequality we
may obtain for positive constants $\lambda$, $\omega_{1}$ and $\omega_{2}$
at our disposal,
\begin{equation}
\begin{alignedat}{1} & \frac{a}{2}\left\Vert C\right\Vert ^{2}+\int_{0}^{t}\left\Vert \nabla C\right\Vert ^{2}ds\\
 & \leq a\left\Vert C_{0}\right\Vert ^{2}+\frac{a}{2}\left\Vert H_{0}\right\Vert ^{2}+\frac{a}{2\lambda}\left\Vert H\right\Vert ^{2}+\frac{\omega_{1}a\psi_{1}}{2}\int_{0}^{t}\oint_{\Gamma}h_{,s}^{2}dAds\\
 & +\frac{bh_{m}d_{5}}{2\omega_{2}\lambda_{1}}+\frac{Ld_{4}}{2}(\gamma_{1}^{-1}+\gamma_{2}^{-1})+\left(\frac{L\gamma_{2}\psi_{1}}{2}+\frac{\psi_{1}K\zeta_{3}}{2}\right)\int_{0}^{t}\oint_{\Gamma}h^{2}dAds\\
 & +\sqrt{\frac{c_{2}}{c_{1}}\int_{0}^{t}\oint_{\Gamma}h^{2}dAds\,\int_{0}^{t}\oint_{\Gamma}|\nabla_{s}h|^{2}dAds}\\
 & +\frac{a\lambda}{2}\left\Vert C\right\Vert ^{2}+\frac{h_{m}b\omega_{2}}{2}\int_{0}^{t}\left\Vert \nabla C\right\Vert ^{2}ds\\
 & +\int_{0}^{t}\left\Vert C\right\Vert ^{2}ds\left(-K+\frac{a}{2\omega_{1}}+\frac{h_{m}b}{\omega_{2}\lambda_{1}^{2}}+\frac{L\gamma_{1}}{2}+\frac{K}{2\zeta_{3}}\right).
\end{alignedat}
\label{eq:proof-26}
\end{equation}
Next, we estimate the $\left\Vert H\right\Vert $ terms using \eqref{eq:h-problem}
and \prettyref{lem:2}, and then pick $\lambda=1/2$ and $\omega_{2}=1/bh_{m}$.
Define the constant $N$ and the data term $d_{6}$ by
\[
\frac{Na}{4}=-K+\frac{a}{2\omega_{1}}+\left(\frac{h_{m}b}{\lambda_{1}}\right)^{2}+\frac{L\gamma_{1}}{2}+\frac{K}{2\zeta_{3}}
\]
and
\[
\begin{alignedat}{1}d_{6} & =a\left\Vert C_{0}\right\Vert ^{2}+\frac{3}{2}\psi_{1}\oint_{\Gamma}h^{2}dA+\frac{bh_{m}d_{5}}{2\omega_{2}\lambda_{1}}+\frac{Ld_{4}}{2}\left(\frac{1}{\gamma_{1}}+\frac{1}{\gamma_{2}}\right)\\
 & +\frac{\omega_{1}a\psi_{1}}{2}\int_{0}^{t}\oint_{\Gamma}h_{,s}^{2}dAds+\sqrt{\frac{c_{2}}{c_{1}}\int_{0}^{t}\oint_{\Gamma}h^{2}dAds\,\int_{0}^{t}\oint_{\Gamma}|\nabla_{s}h|^{2}dAds}\\
 & +\left(\frac{L\gamma_{2}\psi_{1}}{2}+\frac{\psi_{1}K\zeta_{3}}{2}\right)\int_{0}^{t}\oint_{\Gamma}h^{2}dAds.
\end{alignedat}
\]

Then, from \eqref{eq:proof-26} one may derive
\begin{equation}
\frac{a}{4}\left\Vert C\right\Vert ^{2}+\frac{1}{2}\int_{0}^{t}\left\Vert \nabla C\right\Vert ^{2}ds\leq d_{6}+\frac{Na}{4}\int_{0}^{t}\left\Vert C\right\Vert ^{2}ds.\label{eq:proof-27}
\end{equation}
Upon integration of \eqref{eq:proof-27} we may then obtain
\begin{equation}
\int_{0}^{t}\left\Vert C\right\Vert ^{2}ds\leq d_{8}(t)\label{eq:proof-28}
\end{equation}
where $d_{8}$ is the data term
\[
d_{8}=\int_{0}^{t}e^{N(t-s)}d_{6}(s)ds.
\]
Employing \eqref{eq:proof-28} in \eqref{eq:proof-27} one may then
derive the \emph{a priori} bounds
\begin{equation}
\left\Vert C\right\Vert ^{2}\leq Nd_{8}\label{eq:a-priori-1}
\end{equation}
and
\begin{equation}
\int_{0}^{t}\left\Vert \nabla C\right\Vert ^{2}ds\leq\frac{Nad_{8}}{2}.\label{eq:a-priori-2}
\end{equation}

Inequalities \eqref{eq:proof-28} and \eqref{eq:a-priori-1} furnish
the necessary \emph{a priori} bounds for $\left\Vert C\right\Vert $
and we now proceed to derive a similar estimate involving $\left\Vert C\right\Vert _{4}$.

Introduce the function $I(\mathbf{x},t)$ as the solution to 
\begin{equation}
\begin{alignedat}{2}\Delta I & =0 & \qquad & \text{in }\Omega,\\
I & =h^{3}(\mathbf{x},t) &  & \text{on \ensuremath{\Gamma}.}
\end{alignedat}
\end{equation}
Now form the identity
\begin{equation}
\begin{gathered}a\int_{0}^{t}\int_{\Omega}C_{,s}(C^{3}-I)dxds=-b\int_{0}^{t}\int_{\Omega}v_{i}C_{,i}(C^{3}-I)dxds\\
+\int_{0}^{t}\int_{\Omega}\Delta C(C^{3}-I)dxds+L\int_{0}^{t}\int_{\Omega}f(T)(C^{3}-I)dxds\\
-K\int_{0}^{t}\int_{\Omega}C(C^{3}-I)dxds.
\end{gathered}
\label{eq:proof-32}
\end{equation}
Next, denote the five terms in equation \eqref{eq:proof-32} by $J_{1},\ldots,J_{5}.$
We now proceed in a similar manner to that involving the $H$ terms,
employing a weighted arithmetic-geometric mean inequality, Lemmas
\ref{lem:1} and \ref{lem:2}, although we now additionally use Young's
inequality on the term involving $f(T)$. In this way we have
\begin{equation}
J_{1}=\frac{a}{4}(\left\Vert C\right\Vert _{4}^{4}-\left\Vert C_{0}\right\Vert _{4}^{4})-a(I,\, C)+a(I_{0},\, C_{0})+\int_{0}^{t}(I_{,s},\, C)ds.\label{eq:j1}
\end{equation}
Then, 
\begin{equation}
J_{2}=+b\int_{0}^{t}\int_{\Omega}v_{i}C_{,i}Idxds\leq h_{m}^{3}b\sqrt{\int_{0}^{t}\left\Vert \mathbf{v}\right\Vert ^{2}ds\,\int_{0}^{t}\left\Vert \nabla C\right\Vert ^{2}ds}.
\end{equation}
Furthermore,
\begin{equation}
\begin{alignedat}{1}J_{3} & =-\frac{3}{4}\int_{0}^{t}\left\Vert \nabla C^{2}\right\Vert ^{2}ds+\int_{0}^{t}\oint_{\Gamma}h\frac{\partial I}{\partial n}dAds\\
 & \leq-\frac{3}{4}\int_{0}^{t}\left\Vert \nabla C^{2}\right\Vert ^{2}ds+\int_{0}^{t}\sqrt{\oint_{\Gamma}h^{2}dA\,\frac{c_{2}}{c_{1}}\oint_{\Gamma}|\nabla_{s}h|^{2}dA}\, ds.
\end{alignedat}
\end{equation}
For $J_{4}$ and $J_{5}$,
\begin{equation}
\begin{alignedat}{1}J_{4}+J_{5} & \leq\frac{L}{4\varepsilon^{4}}\int_{0}^{t}\left\Vert f\right\Vert _{4}^{4}ds+\frac{3L\varepsilon^{4/3}}{4}\int_{0}^{t}\left\Vert C\right\Vert _{4}^{4}ds\\
 & +\frac{L}{2\varepsilon_{1}}\int_{0}^{t}\left\Vert f\right\Vert ^{2}ds+\left(\frac{\varepsilon_{1}L}{2}+\frac{K}{2\varepsilon_{2}}\right)\int_{0}^{t}\left\Vert I\right\Vert ^{2}ds\\
 & +\frac{K\varepsilon_{2}}{2}\int_{0}^{t}\left\Vert C\right\Vert ^{2}ds-K\int_{0}^{t}\left\Vert C\right\Vert _{4}^{4}ds.
\end{alignedat}
\label{eq:j4-5}
\end{equation}
We group \eqref{eq:j1}-\eqref{eq:j4-5} together in \eqref{eq:proof-32}
to arrive at, after further use of the arithmetic-geometric mean inequality
and \prettyref{lem:2},
\begin{equation}
\begin{alignedat}{1}\frac{a}{4}\left\Vert C\right\Vert _{4}^{4} & +\frac{3}{4}\int_{0}^{t}\left\Vert \nabla C^{2}\right\Vert ^{2}ds+K\int_{0}^{t}\left\Vert C\right\Vert _{4}^{4}ds\\
 & \leq\frac{3a}{4}\left\Vert C_{0}\right\Vert ^{2}+\frac{a\delta_{1}}{2}\left\Vert C\right\Vert ^{2}+\left(\frac{a}{2\delta_{1}}+\frac{a}{2}\right)\psi_{1}\oint_{\Gamma}h^{6}dA\\
 & +\sqrt{\psi_{1}}\int_{0}^{t}\left\Vert C\right\Vert \sqrt{\oint_{\Gamma}|h_{,s}h^{2}|^{2}dA}\, ds\\
 & +\sqrt{\frac{c_{2}}{c_{1}}}\int_{0}^{t}\sqrt{\oint_{\Gamma}h^{2}dA\,\oint_{\Gamma}|\nabla_{s}h^{3}|^{2}dA}\, ds\\
 & +\frac{L}{4\varepsilon^{4}}\int_{0}^{t}\left\Vert f\right\Vert ^{4}ds\\
 & +\frac{L}{2\varepsilon_{1}}\int_{0}^{t}\left\Vert f\right\Vert ^{2}ds+\left(\frac{\varepsilon_{1}L}{2}+\frac{K}{2\varepsilon_{2}}\right)\psi_{1}\int_{0}^{t}\oint_{\Gamma}h^{6}dAds\\
 & +\frac{K\varepsilon_{2}}{2}\int_{0}^{t}\left\Vert C\right\Vert ^{2}ds+h_{m}^{3}b\sqrt{\int_{0}^{t}\left\Vert \nabla C\right\Vert ^{2}ds}\,\sqrt{\frac{d_{5}}{\lambda_{1}}+\frac{2}{\lambda_{1}^{2}}\int_{0}^{t}\left\Vert C\right\Vert ^{2}ds}\\
 & \frac{3L\varepsilon^{4/3}}{4}\int_{0}^{t}\left\Vert C\right\Vert _{4}^{4}ds.
\end{alignedat}
\label{eq:proof-37}
\end{equation}
We now select $\varepsilon=(2K/3L)^{3/4}$. This removes the last
term on the right of \eqref{eq:proof-37}. Now, using the \emph{a
priori} bound on $T$, and the \emph{a priori} estimates \eqref{eq:proof-28}-\eqref{eq:a-priori-2},
we see that what remains on the right hand side of inequality \eqref{eq:proof-37}
may be bounded by data.

Let us denote this term by $\frac{a}{4}d_{9}$. Then
\begin{equation}
\frac{a}{4}\left\Vert C\right\Vert _{4}^{4}+\frac{3}{4}\int_{0}^{t}\left\Vert \nabla C^{2}\right\Vert ^{2}ds+\frac{K}{2}\int_{0}^{t}\left\Vert C\right\Vert _{4}^{4}ds\leq\frac{a}{4}d_{9}(t)
\end{equation}
Inequality \eqref{eq:a-priori-3} thus yields an \emph{a priori} bound
for $\left\Vert C\right\Vert _{4}$, for $\int_{0}^{t}\left\Vert C\right\Vert _{4}^{4}ds$,
and for the term $\int_{0}^{t}\left\Vert \nabla C^{2}\right\Vert ^{2}ds$.
Let us also recollect that using \eqref{eq:proof-17} and \eqref{eq:proof-28}
we have the following \emph{a priori} bound
\begin{equation}
\left\Vert \mathbf{v}\right\Vert ^{2}+\int_{0}^{t}\left\Vert \nabla\mathbf{v}\right\Vert ^{2}ds\leq d_{10},\label{eq:a-priori-3}
\end{equation}
where
\[
d_{10}=d_{5}+\frac{2}{\lambda_{1}}d_{8}.
\]

\section{Continuous dependence on the reaction coefficients}

We now return to the boundary-initial value problem for the difference
equations, \eqref{eq:difference-equations}-\eqref{eq:difference-boundary-conditions}.
In the interests of clarity we recall the specific \emph{a priori}
bounds we require in this section, namely,
\begin{equation}
\begin{alignedat}{1}\int_{0}^{t}\left\Vert \nabla\mathbf{v}\right\Vert ^{2}ds & \leq d_{10}\\
\left\Vert C\right\Vert _{4}^{4} & \leq d_{9}\\
\left\Vert C\right\Vert ^{2} & \leq Nd_{8}
\end{alignedat}
\label{eq:bounds}
\end{equation}

We also require a Sobolev inequality in two-dimensions. This may be
written in the form, see \citet{Payne:1964},
\begin{equation}
\int_{\Omega}|\mathbf{w}|^{4}dx\leq\Omega_{1}\int_{\Omega}|\mathbf{w}|^{2}dx\int_{\Omega}|\nabla\mathbf{w}|^{2}dx.\label{eq:sobolev-inequality}
\end{equation}
where $\Omega_{1}$ is a positive constant. An estimate for $\Omega_{1}$
is given by \citet{Payne:1964}, in Lemma 1, p. 132, as $\Omega_{1}=1/2$.
\begin{thm*}
Let $\chi\equiv(v_{i},p,T,C)$ be a solution to the boundary-initial
value problem \eqref{eq:fundamental}-\eqref{eq:boundary-conditions}
in $\Omega\times(0,\mathscr{T})$ for some $\mathscr{T}<\infty$.
Then the solution $\chi$ depends continuously on the reaction coefficients
$L$ and $K$ explicitly in $L^{2}(\Omega)$ in the sense that the
difference solution $(w_{i},\pi,\theta,\phi)$ given in \eqref{eq:difference-variables}
satisfies the inequality
\[
\left\Vert \mathbf{w}(t)\right\Vert ^{2}+\left\Vert \theta(t)\right\Vert ^{2}+\left\Vert \phi(t)\right\Vert ^{2}\leq f_{1}(t)l^{2}+f_{2}(t)k^{2},\qquad t\in(0,\mathscr{T}),
\]
where $l$, $k$ are defined in \eqref{eq:difference-variables} and
$f_{1}$ and $f_{2}$ are coefficients which depend only on $\Omega$,
$\mathscr{T}$ and the data functions $v_{i}^{0}$, $T_{0}$, $C_{0}$,
$g$ and $h$.\end{thm*}
\begin{proof}
Multiply equation \eqref{eq:difference-equations}$_{1}$ by $w_{i}$
and integrate over $\Omega$ to see that
\[
\begin{alignedat}{1}\frac{d}{dt}\frac{1}{2}\left\Vert \mathbf{w}\right\Vert ^{2} & =-\int_{\Omega}u_{i,j}w_{i}w_{j}dx-\left\Vert \nabla\mathbf{w}\right\Vert ^{2}+(g_{i}\theta,\, w_{i})-(h_{i}\phi,\, w_{i})\\
 & \leq\left\Vert \nabla\mathbf{u}\right\Vert \left\Vert \mathbf{w}\right\Vert _{4}^{2}-\left\Vert \nabla\mathbf{w}\right\Vert ^{2}+\left\Vert \theta\right\Vert \left\Vert \mathbf{w}\right\Vert +\left\Vert \phi\right\Vert \left\Vert \mathbf{w}\right\Vert \\
 & \leq\sqrt{\Omega_{1}}\left\Vert \nabla\mathbf{u}\right\Vert \left\Vert \mathbf{w}\right\Vert \left\Vert \nabla\mathbf{w}\right\Vert -\left\Vert \nabla\mathbf{w}\right\Vert ^{2}+\left\Vert \theta\right\Vert \left\Vert \mathbf{w}\right\Vert +\left\Vert \phi\right\Vert \left\Vert \mathbf{w}\right\Vert ,
\end{alignedat}
\]
where the Cauchy-Schwarz and Sobolev inequalities have been employed.
Now employ the arithmetic-geometric mean inequality for $\beta>0$
to find
\[
\begin{alignedat}{1}\frac{d}{dt}\frac{1}{2}\left\Vert \mathbf{w}\right\Vert ^{2}+\left\Vert \nabla\mathbf{w}\right\Vert ^{2} & \leq\frac{\beta}{2}\left\Vert \mathbf{w}\right\Vert ^{2}\left\Vert \nabla\mathbf{u}\right\Vert ^{2}\\
 & +\frac{\Omega_{1}}{2\beta}\left\Vert \nabla\mathbf{w}\right\Vert ^{2}\\
 & +\frac{\alpha}{2}\left\Vert \theta\right\Vert ^{2}+\frac{\gamma}{2}\left\Vert \phi\right\Vert ^{2}+\left\Vert \mathbf{w}\right\Vert ^{2}\left(\frac{1}{2\alpha}+\frac{1}{2\gamma}\right).
\end{alignedat}
\]
Pick now $\beta=\Omega_{1}$ and employ the Poincaré inequality on
the last term with $\gamma=\alpha$ to obtain
\[
\frac{d}{dt}\frac{1}{2}\left\Vert \mathbf{w}\right\Vert ^{2}+\frac{1}{2}\left\Vert \nabla\mathbf{w}\right\Vert ^{2}\leq\frac{\Omega_{1}}{2}\left\Vert \mathbf{w}\right\Vert ^{2}\left\Vert \nabla\mathbf{u}\right\Vert ^{2}+\frac{\alpha}{2}(\left\Vert \theta\right\Vert ^{2}+\left\Vert \phi\right\Vert ^{2})+\frac{1}{\alpha\lambda_{1}}\left\Vert \nabla\mathbf{w}\right\Vert ^{2}.
\]
Now select $\alpha=4/\lambda_{1}$ and then we find
\begin{equation}
\frac{d}{dt}\left\Vert \mathbf{w}\right\Vert ^{2}+\frac{1}{2}\left\Vert \nabla\mathbf{w}\right\Vert ^{2}\leq\Omega_{1}\left\Vert \nabla\mathbf{u}\right\Vert ^{2}\left\Vert \mathbf{w}\right\Vert ^{2}+\frac{4}{\lambda_{1}}(\left\Vert \theta\right\Vert ^{2}+\left\Vert \phi\right\Vert ^{2}).\label{eq:proof-42}
\end{equation}

Upon multiplying equation \eqref{eq:difference-equations}$_{3}$
by $\theta$ and integrating over $\Omega$ we derive
\begin{equation}
\begin{alignedat}{1}\frac{d}{dt}\frac{1}{2}\left\Vert \theta\right\Vert ^{2} & =-\left\Vert \nabla\theta\right\Vert ^{2}+\int_{\Omega}w_{i}T_{1}\theta_{,i}dx\\
 & \leq-\left\Vert \nabla\theta\right\Vert ^{2}+T_{m}\left\Vert \mathbf{w}\right\Vert \left\Vert \nabla\theta\right\Vert \\
 & \leq\frac{T_{m}^{2}}{4}\left\Vert \mathbf{w}\right\Vert ^{2}.
\end{alignedat}
\label{eq:proof-43}
\end{equation}

The next stage requires us to multiply equation \eqref{eq:difference-equations}$_{4}$
by $\phi$ and integrate over $\Omega$ to find
\begin{equation}
\begin{alignedat}{1}\frac{a}{2}\frac{d}{dt}\left\Vert \phi\right\Vert ^{2} & =b\int_{\Omega}w_{i}C_{1}\phi_{,i}dx-\left\Vert \nabla\phi\right\Vert ^{2}-K_{1}\left\Vert \phi\right\Vert ^{2}-k(C_{2},\,\phi)\\
 & +L_{1}(f(T^{1})-f(T^{2}),\,\phi)+l(f(T^{2}),\,\phi).
\end{alignedat}
\label{eq:proof-44}
\end{equation}
Using Lagrange's theorem we know $f(T^{1})-f(T^{2})=\theta f'(\xi)$
for some $\xi\in(T_{1},T_{2})$. Then since $T_{m}$ is a bound for
$|T|$ and $f\in C^{1}$, we know $|f'(\xi)|\leq d_{1}$, $|f(T^{2})|\leq d_{2}$
for data terms $d_{1}$ and $d_{2}$. Thus, the last two terms of
\eqref{eq:proof-44} may be bounded by
\begin{equation}
L_{1}d_{1}\left(\frac{\left\Vert \theta\right\Vert ^{2}}{2\alpha}+\frac{\alpha}{2}\left\Vert \phi\right\Vert ^{2}\right)+\frac{l^{2}\left\Vert f(T_{2})\right\Vert ^{2}}{2\beta}+\frac{\beta}{2}\left\Vert \phi\right\Vert ^{2},\label{eq:proof-45}
\end{equation}
for $\alpha,\beta>0$ to be selected. Likewise for $\gamma>0$,
\begin{equation}
-k(C_{2},\,\phi)\leq k^{2}\frac{\left\Vert C_{2}\right\Vert ^{2}}{2\gamma}+\frac{\gamma}{2}\left\Vert \phi\right\Vert ^{2}.
\end{equation}

For the cubic term we have
\begin{equation}
\begin{alignedat}{1}b\int_{\Omega}w_{i}C_{1}\phi_{,i}dx & \leq\frac{b}{2\zeta}\int_{\Omega}|\mathbf{w}|^{2}C_{1}^{2}dx+\frac{\zeta b}{2}\left\Vert \nabla\phi\right\Vert ^{2}\\
 & \leq\frac{b}{2\zeta}\left\Vert \mathbf{w}\right\Vert _{4}^{2}\left\Vert C_{1}\right\Vert _{4}^{2}+\frac{\zeta b}{2}\left\Vert \nabla\phi\right\Vert ^{2}\\
 & \leq\frac{b\sqrt{\Omega_{1}}}{2\zeta}\left\Vert C_{1}\right\Vert _{4}^{2}\left\Vert \mathbf{w}\right\Vert \left\Vert \nabla\mathbf{w}\right\Vert +\frac{\zeta b}{2}\left\Vert \nabla\phi\right\Vert ^{2},\\
 & \leq\frac{b\Omega_{1}}{4\zeta\mu}\left\Vert C_{1}\right\Vert _{4}^{4}\left\Vert \mathbf{w}\right\Vert ^{2}+\frac{b\mu}{4\zeta}\left\Vert \nabla\mathbf{w}\right\Vert ^{2}+\frac{\zeta b}{2}\left\Vert \nabla\phi\right\Vert ^{2},
\end{alignedat}
\label{eq:proof-47}
\end{equation}
for $\zeta,\mu>0$ to be chosen.

We employ \eqref{eq:proof-45}-\eqref{eq:proof-47} in \eqref{eq:proof-44},
pick $\beta=\gamma=K_{1}/2$, $\alpha=K_{1}/L_{1}d_{1}$ and pick
$\zeta=2/b$. Then we obtain
\begin{equation}
\begin{alignedat}{1}\frac{a}{2}\frac{d}{dt}\left\Vert \phi\right\Vert ^{2} & \leq\frac{L_{1}^{2}d_{1}^{2}}{2K_{1}}\left\Vert \theta\right\Vert ^{2}+\frac{\left\Vert f(T^{2})\right\Vert ^{2}}{K_{1}}l^{2}+\frac{\left\Vert C_{2}\right\Vert ^{2}}{K_{1}}k^{2}\\
 & +\frac{\Omega_{1}b^{2}}{8\mu}\left\Vert C_{1}\right\Vert _{4}^{4}\left\Vert \mathbf{w}\right\Vert ^{2}+\frac{b^{2}\mu}{8}\left\Vert \nabla\mathbf{w}\right\Vert ^{2}.
\end{alignedat}
\label{eq:proof-48}
\end{equation}

Thus, from \eqref{eq:proof-42}, \eqref{eq:proof-43} and \eqref{eq:proof-48}
we may obtain
\begin{equation}
\begin{alignedat}{1} & \frac{d}{dt}(\left\Vert \mathbf{w}(t)\right\Vert ^{2}+\left\Vert \theta(t)\right\Vert ^{2}+\left\Vert \phi(t)\right\Vert ^{2})+\frac{1}{2}\left\Vert \nabla\mathbf{w}\right\Vert ^{2}\\
 & \qquad\leq\left(\frac{T_{m}^{2}}{2}+\frac{\Omega_{1}b^{2}}{4a\mu}\left\Vert C_{1}\right\Vert _{4}^{4}+\Omega_{1}\left\Vert \nabla\mathbf{u}\right\Vert ^{2}\right)\\
 & \qquad\qquad+\frac{b^{2}\mu}{4a}\left\Vert \nabla\mathbf{w}\right\Vert ^{2}+\left(\frac{4}{\lambda_{1}}+\frac{L_{1}^{2}d_{1}^{2}}{aK_{1}}\right)\left\Vert \theta\right\Vert ^{2}+\frac{4}{\lambda_{1}}\left\Vert \phi\right\Vert ^{2},\\
 & \qquad\qquad+\frac{2}{aK_{1}}\left\Vert f(T^{2})\right\Vert ^{2}l^{2}+\frac{2\left\Vert C_{2}\right\Vert ^{2}}{aK_{1}}k^{2}.
\end{alignedat}
\label{eq:proof-49}
\end{equation}

Use the bounds \eqref{eq:bounds} on $\left\Vert C_{1}\right\Vert _{4}$
and $\left\Vert C_{2}\right\Vert $, choose $\mu=2a/b^{2}$, and set
\[
M=\max\left\{ \frac{\Omega_{1}b^{4}d_{9}}{8a^{2}}+\frac{T_{m}^{2}}{2},\,\frac{4}{\lambda_{1}}+\frac{L_{1}^{2}d_{1}^{2}}{aK_{1}}\right\} .
\]

The term $\left\Vert f(T^{2})\right\Vert $ is \emph{a priori} bounded
and so we set
\[
\alpha_{1}=\frac{2}{aK_{1}}\left\Vert f(T_{2})\right\Vert ^{2},\qquad\alpha_{2}=\frac{2}{aK_{1}}\left\Vert C_{2}\right\Vert ^{2}.
\]

Then, with $F(t)=\left\Vert \mathbf{w}(t)\right\Vert ^{2}+\left\Vert \theta(t)\right\Vert ^{2}+\left\Vert \phi(t)\right\Vert ^{2}$
we obtain from \eqref{eq:proof-49}
\begin{equation}
F'-(\Omega_{1}\left\Vert \nabla\mathbf{u}\right\Vert ^{2}+M)F\leq\alpha_{1}l^{2}+\alpha_{2}k^{2}.\label{eq:differential-equation}
\end{equation}
We integrate \eqref{eq:differential-equation} via an integrating
factor and we use \eqref{eq:bounds}$_{1}$ to now see that
\begin{equation}
\left\Vert \mathbf{w}(t)\right\Vert ^{2}+\left\Vert \theta(t)\right\Vert ^{2}+\left\Vert \phi(t)\right\Vert ^{2}\leq R(t)(\alpha_{1}l^{2}+\alpha_{2}k^{2}),\label{eq:energy-inequality}
\end{equation}
where $R(t)$ is the data term
\begin{equation}
R(t)=\int_{0}^{t}\exp\left[M(t-s)+\Omega_{1}\int_{s}^{t}d_{10}(y)dy\right]ds.
\end{equation}

Inequality \eqref{eq:energy-inequality} is the sought after result
and establishes continuous dependence on the reaction coefficients.\end{proof}
\begin{rem*}
We have established continuous dependence when $\Omega$ is a bounded
domain in $\mathbb{R}^{2}$. The proof relies on the Sobolev inequality
\eqref{eq:sobolev-inequality}. This inequality does not hold when
$\Omega\subset\mathbb{R}^{3}$ and instead one must employ an alternative
Sobolev inequality of form
\[
\left\Vert \mathbf{w}\right\Vert _{4}^{4}\leq\tilde{\Omega}_{1}\left\Vert \mathbf{w}\right\Vert \left\Vert \nabla\mathbf{w}\right\Vert ^{3}.
\]
In this case we cannot establish continuous dependence by methods
analogous to those employed here. Evidently we may proceed along similar
lines as \citet{PayneStraughan:1989} for the three-dimensional problem
but only for $t$ restricted. We do not give details of such a local
result since we are interested in a truly global \emph{a priori} result.
\end{rem*}
\bibliographystyle{unsrtnat}
\bibliography{chemotaxis}

\end{document}